# Numerical solution of the time fractional nonlinear Fisher-KPP diffusion-reaction equation using the local domain boundary element method

Theodore V. Gortsas

*Department of Mechanical Engineering and Aeronautics, University of Patras, Patras, Greece*

**ABSTRACT**

The Fisher-KPP partial differential equation has been employed in science to model various biological, chemical, and thermal phenomena. Time fractional extensions of Fisher's equation have also appeared in the literature, aiming to model systems with memory. The solution of the time fractional Fisher-KPP equation is challenging due to the interplay between the nonlinearity and the nonlocality imposed by the fractional derivatives. An accurate method that for the solution of time fractional diffusion problems is the Boundary Element Method (BEM). The conventional BEM has a high computational cost and memory requirements since it leads to dense coefficient matrices. For nonlinear transient problems, its efficiency is further reduced due to the appearance of volume integrals. In the present work an extension of the recently proposed Local Domain Boundary Element Method (LD-BEM) is presented for the solution of nonlinear time fractional Fisher-KPP problems. The implemented numerical method is used to examine various two-dimensional problems related to the Fisher-KPP equation using different definitions of the fractional derivative.

## 1. Introduction

Fractional calculus [1,2] has emerged as an important mathematical tool in the physical sciences, extending the modelling power of partial differential equations. Its applications span diverse fields, such as biology [3-6], as well as material science [7-9] and engineering [10,11]. Many of the modelled systems in these fields exhibit to a certain level non-locality, memory, spatial heterogeneity as well as anomalous diffusion [12]. The study of fractional differential equations has evolved along two different directions, one being the understanding of the physical basis that justifies their use as models [13,14] and the other being the development of the numerical tools required for their solution [15]. One of the challenges when modelling various phenomena with fractional differential equations, is the choice of the fractional derivatives, since the definition of fractional derivatives is not unique, and many different types of derivatives have appeared in the literature [16,17].

The nonlinear Fisher-KPP partial differential equation [18, 19] has been largely used to model various biological phenomena [18-28], as well as various thermal and chemical processes [29-31]. Time fractional extensions of Fisher's equation have also been presented, aiming to model systems with memory. Such time fractional extensions have been discussed in the recent works of Alquran et al. [32] and Angstmann and Henry [33].The solution of the time fractional Fisher-KPP equation is challenging due to the interplay between the nonlinearity and the nonlocality imposed by the fractional derivatives. The analytical, semi-analytical, and computational approaches that have been developed, mainly concern one-dimensional problems [34-41]. Only a limited number of studies, have addressed two-dimensional time-fractional Fisher-KPP boundary value problems [42-44].

An accurate method that has been proven effective for the solution of time fractional diffusion problems is the Boundary Element Method (BEM). In Carrer et al. [45] the BEM has been applied for the solution of the time fractional wave equation using the Caputo fractional derivative, while in Dehghan and Safarpoor [46] it has been applied to the solution of some nonlinear two-dimensional time fractional partial differential equations using both the Caputo and the Riemann-Liouville fractional derivatives. One of the challenges of BEM is the computation of the large number of volume integrals which appear in the formulation of the numerical method, since a fundamental solution does not exist for the nonlinear differential operator. The computation the volume integrals is usually achieved through direct integration by using a volume mesh, or by using techniques such as the dual reciprocity method [46-48] and the radial integration method [49-51], which efficiently transform the domain integrals appearing in the BEM formulations into equivalent boundary integrals, thus allowing the boundary-only discretization of the standard BEM. Another challenge is the storage and the solution of the final linear system of equations, as the resulting matrix is dense and fully populated. For this reason, techniques based on the Fast Multipole Method [52] and Hierarchical Matrices [53-55] have been proposed which in combination with iterative solvers aim to reduce the computational cost and the memory requirements of the method. The efficient implementation of such methods for nonlinear problems, is a highly complex task.

In the present work, the recently proposed Local Domain Boundary Element Method (LD-BEM) [56] is extended for the solution of nonlinear time fractional Fisher-KPP problems. The LD-BEM employs a poly-region discretization, leading to the solution of a linear system of equations with a sparse system matrix. Moreover, the number of volume integrals that must be computed is significantly reduced, while the uknown degrees of freedom are limited to those on the boundary of each subregion. The implemented numerical method is applied to several problems with different definitions of the fractional derivative, such as the Caputo fractional derivative [16] as well as a fractal fractional derivative recently proposed in [17]. Except the LD-BEM, for the sake of comparison, the meshless Fragile Points Method (FPM) [57] has been implemented for the numerical solution of nonlinear time fractional Fisher-KPP problems.

The structure of the paper is the following: In Section 2 the mathematical setup of the boundary value problem is outlined, along with the definitions of the fractional derivatives employed in this work. Furthermore, the integral representation of the solution is presented. In Section 3 the steps for the numerical implementation of LD-BEM for the solution of the Fisher-KPP equation are described. In Section 4 the implemented numerical method is used for the solution of a series of numerical examples. Finally, in Section 5 the main conclusions of this work are presented.

## 2. Time fractional Fisher-KPP equation and integral representation of the solution

### 2.1 *Time fractional Fisher-KPP diffusion equation*

The time fractional Fisher-KPP equation can be expressed as follows

$$\nabla^2 \varphi(\mathbf{r},t) = \frac{1}{\rho}\left( \frac{\partial^{a,\beta} \varphi(\mathbf{r},t)}{\partial t^{a,\beta}} - N(\{\mathbf{m}\},\varphi(\mathbf{x},t)) - f(\mathbf{r},t) \right), \quad (\mathbf{r},t) \in V \times (0,t_D] \quad (1)$$

where $\varphi(\mathbf{r},t)$ denotes the concentration function, $\rho \in \mathbb{R}^+$ is the diffusion coefficient, $\frac{\partial^{a,\beta} \varphi(\mathbf{r},t)}{\partial t^{a,\beta}}$ denotes a general fractional or fractal fractional time derivative with $a, \beta \in (0,1)$ and $V \subseteq \mathbb{R}^2$ is a domain with boundary $\partial V$. The function $f : V \times [0,t] \to \mathbb{R}$ is a predefined function, while $N(\mathbf{m},\varphi)$ denotes a nonlinear term with $\mathbf{m}=[b,c,d]^T \in \mathbb{R}^3$ being a constant coefficient vector. The nonlinear term $N(\mathbf{m},\varphi)$ obtains the following form,

$$N(\mathbf{m},\varphi(\mathbf{r},t)) = \varphi(\mathbf{r},t)\left(c - d\varphi^b(\mathbf{r},t)\right) \qquad (2)$$

The boundary and initial conditions are given as follows,

$$\begin{aligned}
\varphi(\mathbf{r},t) &= \bar{\varphi}(\mathbf{r},t), \quad \mathbf{r} \in \partial V_\varphi \\
q(\mathbf{r},t) &= \mathbf{n}(\mathbf{r}) \cdot \nabla \varphi(\mathbf{r},t) = \bar{q}(\mathbf{r},t), \quad \mathbf{r} \in \partial V_q \\
\varphi(\mathbf{r},0) &= \varphi_0(\mathbf{r}), \quad \mathbf{r} \in V
\end{aligned} \qquad (3)$$

where the parts of the boundary $\partial V_q$ and $\partial V_\varphi$ are considered disjoint, $\mathbf{n}(\mathbf{r})$ is the unit normal vector at $\mathbf{r}$ and $\bar{\varphi}(\mathbf{r},t)$, $\bar{q}(\mathbf{r},t)$, $\varphi_0(\mathbf{r})$ are prescribed functions.

2.2 *Definitions of fractional time derivative*

In this work multiple expressions for the time fractional derivative appearing in Eq. (1) are examined. The first one is the Caputo fractional derivative, which is defined as [16, 35],

$$\frac{\partial \varphi^a(\mathbf{r},t)}{\partial t^a} = {}^C D_t^a(\varphi) = \frac{1}{\Gamma(1-\alpha)} \int_0^t K_C(t-\tau,\alpha)\dot{\varphi}(\mathbf{r},\tau)d\tau, \quad a \in (0,1) \qquad (4)$$

where $K_c(t-\tau,a) = (t-\tau)^{-\alpha}$ is the kernel of the integral operator of Eq. (4) and $\Gamma(a)$ denotes the gamma function defined as $\Gamma(\alpha) = \int_0^\infty \tau^{\alpha-1}e^{-t}dt$ for $a > 0$.

The second expression examined for the time derivative of the function appearing in Eq. (1) is the fractal fractional derivative of order $\alpha,\beta$ in the Riemann-Liouville sense, recently introduced in [17], which aims to combine a local arbitrary differentiation order with the nonlocal definition of the derivative. This expression is the following,

$$\frac{\partial \varphi^{a,\beta}(\mathbf{r},t)}{\partial t^{a,\beta}} = {}^{FFP}D_t^{a,\beta}(\varphi) = \frac{1}{\Gamma(1-\alpha)} \frac{d}{dt^\beta} \int_0^t K_C(t-\tau,a)\varphi(\mathbf{r},\tau)d\tau, \quad a,\beta \in (0,1) \qquad (5)$$

where $\frac{d}{dt^\beta}(\varphi(\mathbf{r},t)) = \lim_{t_0 \to t} \frac{\varphi(\mathbf{r},t_0) - \varphi(\mathbf{r},t)}{(t_0)^\beta - (t)^\beta}$ and the kernel of the integral operator is the same as the one appearing in the definition of the Caputo fractional derivative. Furthermore, for $\beta = 1$ Eq. (5) leads to the Riemann-Liouville fractional derivative defined as,

$$\frac{\partial \varphi^{a,1}(\mathbf{r},t)}{\partial t^{a,1}} = {}^{RL}D_t^a(\varphi) = \frac{1}{\Gamma(1-\alpha)} \frac{d}{dt} \int_0^t K_C(t-\tau, a)\varphi(\mathbf{r},\tau)d\tau, \quad a \in (0,1) \tag{6}$$

A useful relation between the Caputo and the Riemann-Liouville fractional time derivatives, which has been used in [17] is the following,

$$^{RL}D_t^a(\varphi) = {}^{C}D_t^a(\varphi) + \frac{t^{-a}\varphi(\mathbf{r},0)}{\Gamma(1-\alpha)}, \quad a \in (0,1) \tag{7}$$

2.3 *Integral representation of the solution*

To develop an integral formulation for a transient problem with the BEM two different approaches are commonly employed. The first approach is to use a fundamental solution for the time-space operator, which leads to a boundary only formulation [58] if there are no volume sources. The second approach, which is the one followed in this work, is to use the fundamental solution of the related static problem [50, 51], an approach which leads to the appearance of volume integrals. Moreover, to develop an integral formulation for nonlinear problems, it is common to split the differential operator to linear and nonlinear parts [59] and use the fundamental solution of the linear part.

Using the fundamental solution $G(\mathbf{r},\mathbf{r}')$ of the Laplace equation [60], and considering a time discretization of the interval $[0, t_D]$, an integral representation of the solution of Eq. (1) can be derived at $t_{n+1}$, given as follows,

$$\begin{aligned}
&c(\mathbf{r})\varphi_{n+1}(\mathbf{r}) + \int_{\partial V} Q(\mathbf{r},\mathbf{r}')\varphi_{n+1}(\mathbf{r}')dS = \int_{\partial V} G(\mathbf{r},\mathbf{r}')q_{n+1}(\mathbf{r}')dS - \\
&\frac{1}{\rho}\int_V G(\mathbf{r},\mathbf{r}')\left(\frac{\partial^a \varphi}{\partial t^a}(\mathbf{r}')\right)_{n+1} dV + \frac{1}{\rho}\int_V G(\mathbf{r},\mathbf{r}')N(\{\mathbf{m}\},\varphi_{n+1}(\mathbf{r}'))V + \\
&\frac{1}{\rho}\int_V G(\mathbf{r},\mathbf{r}')f_{n+1}(\mathbf{r}')dV, \quad n = 0\ldots N_t - 1
\end{aligned} \tag{8}$$

where the value of the jump coefficient $c(\mathbf{r})$ depends on location of the collocation point $\mathbf{r}$. For collocation points on smooth parts of the boundaries $c(\mathbf{r}) = 0.5$, while for interior points $c(\mathbf{r}) = 1$. The fundamental solution $G(\mathbf{r},\mathbf{r}')$ and its derivative $Q(\mathbf{r},\mathbf{r}')$ are given as,

$$G(\mathbf{r},\mathbf{r}') = -\frac{1}{2\pi} \ln(|\mathbf{r} - \mathbf{r}'|)$$

$$Q(\mathbf{r},\mathbf{r}') = \mathbf{n}' \cdot \nabla_{\mathbf{r}'} G(\mathbf{r},\mathbf{r}')$$

(9)

with |$\mathbf{r}$ - $\mathbf{r}'$| being the distance between the points $\mathbf{r}$ and $\mathbf{r}'$.

Using directly Eq. (8) for the solution of the problem leads to a quadratic computational cost since the fundamental solution has global support. For this reason, the domain $V$ is partitioned in conformal, non-overlapping subregions as described in [56] and the integral equation is applied separately for each subregion. On the interfaces between the subregions the following conditions are applied

$$\varphi_1 - \varphi_2 = 0, \quad q_1 + q_2 = 0 \tag{10}$$

where the indices indicate the two sides of the interface with the normal vectors $\mathbf{n}_2 = -\mathbf{n}_1$. The above conditions ensure the continuity of the concentration function and the flux equilibrium.

## 3. Numerical implementation of the LD-BEM

3.1 *Discretization of time fractional derivatives*

The fractional time derivative in the Caputo sense, considering the time discretization can be expressed at $t_{n+1}$ as,

$$^{C}D^{a}_{t_{n+1}}(\varphi) = \frac{1}{\Gamma(1-\alpha)} \sum_{k=0}^{n} \int_{t_k}^{t_{k+1}} (t_{n+1} - \tau)^{-a} \dot{\varphi}(\tau) d\tau \tag{11}$$

Separating the time interval $[t_n, t_{n+1}]$ in the summation, assuming that the time derivative is constant within each interval the following equation is obtained,

$$^{C}D_{t_{n+1}}^{a}(\varphi) = \frac{1}{\Gamma(1-\alpha)}\left(\left[-\frac{(t_{n+1}-\tau)^{1-a}}{1-a}\right]_{t_n}^{t_{n+1}}\dot{\varphi}_{n+1} + \sum_{k=0}^{n-1}\left[-\frac{(t_{n+1}-\tau)^{1-a}}{1-a}\right]_{t_k}^{t_{k+1}}\dot{\varphi}_{k+1}\right) \quad (12)$$

The first term in the above expression can be expressed as,

$$\left[-\frac{(t_{n+1}-\tau)^{1-a}}{(1-a)}\right]_{t_n}^{t_{n+1}}\dot{\varphi}_{n+1} = \frac{\Delta t}{(1-a)\Delta t^a}\dot{\varphi}_{n+1} \quad (13)$$

while the second term can be rewritten after some algebra as follows,

$$\sum_{k=0}^{n-1}\left[-\frac{(t_{n+1}-\tau)^{1-a}}{(1-a)}\right]_{t_k}^{t_{k+1}}\dot{\varphi}_{k+1} = \sum_{k=0}^{n-1}\frac{\Delta t}{(1-a)\Delta t^a}\left[(n+1-k)^{1-a}-(n-k)^{1-a}\right]\dot{\varphi}_{k+1} \quad (14)$$

Combining Eqs. (12), (13) and (14) and utilizing the property of the gamma function $\Gamma(2-a) = \Gamma(1-a)(1-a)$ the following approximation of the fractional derivative can be obtained,

$$^{C}D_{t_{n+1}}^{a}(\varphi) = \frac{\Delta t}{\Delta t^a \Gamma(2-a)}\left(\dot{\varphi}_{n+1} + \sum_{k=0}^{n-1}B_{n,k,a}\dot{\varphi}_{k+1}\right) \quad (15)$$

where $B_{n,k,a} = (n+1-k)^{1-a} - (n-k)^{1-a}$. Introducing the finite difference approximation of the derivative $\dot{\varphi}_{n+1} = \frac{\varphi_{n+1}-\varphi_n}{\Delta t}$ the following equation is derived,

$$^{C}D_{t_{n+1}}^{a}(\varphi) = c_a\left(\varphi_{n+1} - \varphi_n + P_{n,a}\right) \quad (16)$$

where $c_a = \frac{1}{\Delta t^a \Gamma(2-a)}$, and $P_{n,a} = \sum_{k=0}^{n-1}B_{n,k,a}\left(\varphi_{k+1} - \varphi_k\right)$.

The fractal fractional time derivative in the Riemann-Liouville sense, considering the properties of the fractal derivative and the relation between the Riemann-Liouville and Caputo fractional derivatives can be written at $t_{n+1}$ as,

$$^{FFP}D_{t_{n+1}}^{a,\beta}(\varphi) = \frac{1}{\beta\,\Gamma(1-\alpha)}\frac{1}{(t_{n+1})^{\beta-1}}\left(\sum_{k=0}^{n}\int_{t_k}^{t_{k+1}}(t_{n+1}-\tau)^{-a}\dot{\varphi}(\tau)d\tau + \frac{t_{n+1}^{-a}\varphi(0)}{\Gamma(1-\alpha)}\right) \quad (17)$$

Separating the time interval $[t_n, t_{n+1}]$ in the summation the following equation is obtained,

$$^{FFP}D_{t_{n+1}}^{a,\beta}(\varphi) = \frac{1}{\beta\,\Gamma(1-\alpha)}\frac{1}{(t_{n+1})^{\beta-1}}\left(\frac{t_{n+1}^{-a}\varphi(0)}{\Gamma(1-\alpha)} + \dot{\varphi}_{n+1}\int_{t_n}^{t_{n+1}}(t_{n+1}-\tau)^{-a}\,d\tau \right. \\ \left. + \sum_{k=0}^{n-1}\dot{\varphi}_{k+1}\int_{t_k}^{t_{k+1}}(t_{n+1}-\tau)^{-a}\,d\tau\right) \tag{18}$$

The second and third terms in the summation are the same as the ones appearing in the discretization of the Caputo fraction derivative, and they can be replaced by Eq. (13) and Eq. (14) respectively,

$$^{FFP}D_{t_{n+1}}^{a,\beta}(\varphi) = \frac{1}{\beta\,\Gamma(1-\alpha)}\frac{\Delta t}{\Delta t^{\beta}(n+1)^{\beta-1}}\left(\frac{\varphi(0)}{\Delta t^a(n+1)^a\Gamma(1-\alpha)} + \right. \\ \left. \frac{\Delta t}{(1-a)\Delta t^a}\dot{\varphi}_{n+1} + \sum_{k=0}^{n-1}\frac{\Delta t}{(1-a)\Delta t^a}B_{n,k,a}\dot{\varphi}_{k+1}\right) \tag{19}$$

Introducing the finite difference approximation of the derivative $\dot{\varphi}_{n+1} = \frac{\varphi_{n+1} - \varphi_n}{\Delta t}$ the following equation is obtained,

$$^{FFP}D_{t_{n+1}}^{a,\beta}(\varphi) = c_{a,\beta,n}\left(\varphi_{n+1} - \varphi_n + P_{n,a,\beta}\right) \tag{20}$$

Where the constant $c_{a,\beta,n}$ is given as $c_{a,\beta,n} = \frac{1}{\Gamma(2-\alpha)}\frac{\Delta t}{\beta\,\Delta t^{\alpha+\beta}(n+1)^{\beta-1}}$, and the term $P_{n,a,\beta}$ is given as $P_{n,a,\beta} = \sum_{k=0}^{n-1}B_{n,k,a}\left(\varphi_{k+1} - \varphi_k\right) + \frac{\varphi(0)(1-a)}{(n+1)^a\Gamma(1-\alpha)}$.

### 3.2 Discrete equations and numerical solution

In this section, the discrete system of equations is presented for a single subdomain. After these equations been derived, an assembly procedure needs to be implemented to obtain the final sparse linear system of equations. Considering a subregion to be of quadrilateral shape, we assume without loss of generality that the boundary edges of each subregion are discretized with linear discontinuous boundary elements [61]. Furthermore, we consider each subregion to be a linear discontinuous cell element.

Collocating Eq. (8) for the nodal points related to the boundary discretization, and for the nodal points related to the volume discretization we obtain the following equations,

$$\frac{1}{2}\varphi_{n+1}^k + \sum_{l=1}^{N_E^{(D)}} \sum_{b=1}^{2} H^{klb} \varphi_{n+1}^b = \sum_{l=1}^{N_E^{(D)}} \sum_{b=1}^{2} G^{klb} q_{n+1}^b - \frac{1}{\rho} \sum_{l=1}^{N_E^{(D)}} C^{ki} \left(\frac{\partial^a \varphi}{\partial t^a}\right)_{n+1}^i +$$

$$\frac{1}{\rho} \sum_{l=1}^{N_E^{(D)}} C^{ki} N^i\left(\mathbf{m}, \varphi_{n+1}^i\right) + \frac{1}{\rho} \sum_{l=1}^{N_E^{(D)}} C^{ki} f_{n+1}^i \qquad (21)$$

$$\varphi_{n+1}^j + \sum_{l=1}^{N_E^{(D)}} \sum_{b=1}^{2} H^{jlb} \varphi_{n+1}^b = \sum_{l=1}^{N_E^{(D)}} \sum_{b=1}^{2} G^{jlb} q_{n+1}^b - \frac{1}{\rho} \sum_{l=1}^{N_E^{(D)}} C^{ji} \left(\frac{\partial^a \varphi}{\partial t^a}\right)_{n+1}^i +$$

$$\frac{1}{\rho} \sum_{l=1}^{N_E^{(D)}} C^{ji} N^i\left(\mathbf{m}, \varphi_{n+1}^i\right) + \frac{1}{\rho} \sum_{l=1}^{N_E^{(D)}} C^{ji} f_{n+1}^i \qquad (22)$$

The quantities appearing in Eq. (21) can be calculated through the following integrations,

$$H^{klb} = \int_{-1}^{1} Q\left(\mathbf{z}_{(\partial D)}^k, \mathbf{r}'(\xi)\right) N^b(\xi) |J^l(\xi)| d\xi$$

$$G^{klb} = \int_{-1}^{1} G\left(\mathbf{z}_{(\partial D)}^k, \mathbf{r}'(\xi)\right) N^b(\xi) |J^l(\xi)| d\xi \qquad (23)$$

$$C^{ki} = \int_{-1}^{1}\int_{-1}^{1} G\left(\mathbf{z}_{(\partial D)}^k, \mathbf{r}'(\xi, n)\right) N^i(\xi, n) |J(\xi)| d\xi dn$$

where the symbol $\mathbf{z}_{(\partial D)}^k$ denotes a boundary collocation point. The quantities appearing in Eq. (22) can be expressed in a similar way using the collocation points $\mathbf{z}_{(D)}^i$ of the discontinuous cell element. When a collocation point corresponds to a nodal point in the element of integration the integrals of Eq. (23) exhibit singular behavior. For the computation of these integrals the semi-analytical integration procedures of Guiggiani and Gigante [62] are utilized.

As it was shown in Eqs. (16) and (20) the general form for the discrete fractional time derivative can be expressed as follows,

$$\left(\frac{\partial^a \varphi}{\partial t^a}\right)_{n+1} = c_a \left(\varphi_{n+1} - \varphi_n + P_{n,a}\right) \tag{24}$$

Substituting the previous equation to Eqs. (21) and (22), and rewriting the equations in matrix form we obtain,

$$\begin{bmatrix} \mathbf{H}^{(BB)} & \mathbf{H}^{(BI)} \\ \mathbf{H}^{(IB)} & \mathbf{H}^{(II)} \end{bmatrix} \begin{Bmatrix} \boldsymbol{\varphi}_{n+1}^{(B)} \\ \boldsymbol{\varphi}_{n+1}^{(I)} \end{Bmatrix} = \begin{bmatrix} \mathbf{G}^{(BB)} & \mathbf{G}^{(BI)} \\ \mathbf{G}^{(IB)} & \mathbf{G}^{(II)} \end{bmatrix} \begin{Bmatrix} \mathbf{q}_{n+1}^{(B)} \\ \mathbf{q}_{n+1}^{(I)} \end{Bmatrix}$$

$$-\frac{1}{\rho}\begin{bmatrix} \mathbf{C}^{(BD)} \\ \mathbf{C}^{(ID)} \end{bmatrix} \left( c_a \left( \{\boldsymbol{\varphi}_{n+1}^{(D)}\} - \{\boldsymbol{\varphi}_n^{(D)}\} + \{\mathbf{P}_{n,a}^{(D)}\} \right) - \{\mathbf{N}(\mathbf{m}, \boldsymbol{\varphi}_{n+1}^{(D)})\} - \{\mathbf{f}_{n+1}^{(D)}\} \right) \tag{25}$$

$$\{\boldsymbol{\varphi}_{n+1}^{(D)}\} + \begin{bmatrix} \mathbf{H}^{(DB)} & \mathbf{H}^{(DI)} \end{bmatrix} \begin{Bmatrix} \boldsymbol{\varphi}_{n+1}^{(B)} \\ \boldsymbol{\varphi}_{n+1}^{(I)} \end{Bmatrix} = \begin{bmatrix} \mathbf{G}^{(DB)} & \mathbf{G}^{(DI)} \end{bmatrix} \begin{Bmatrix} \mathbf{q}_{n+1}^{(B)} \\ \mathbf{q}_{n+1}^{(I)} \end{Bmatrix}$$

$$-\frac{1}{\rho}\begin{bmatrix} \mathbf{C}^{(DD)} \end{bmatrix} \left( c_a \left( \{\boldsymbol{\varphi}_{n+1}^{(D)}\} - \{\boldsymbol{\varphi}_n^{(D)}\} + \{\mathbf{P}_{n,a}^{(D)}\} \right) - \{\mathbf{N}(\mathbf{m}, \boldsymbol{\varphi}_{n+1}^{(D)})\} - \{\mathbf{f}_{n+1}^{(D)}\} \right) \tag{26}$$

The partitioned matrices [**H**] and [**G**] which appear in Eq. (25) contain the boundary integrals of Eq. (21). Since linear discontinuous elements have been used for the discretization of the boundary edges their dimensions are $8 \times 8$. The partitioned matrices [**H**] and [**G**] appearing in Eq. (26) are rectangular matrices of dimension $4 \times 8$ and contain the boundary integrals of Eq. (22). The matrix [**C**] contains the volume integrals of each equation and is of size $8 \times 4$ for Eq. (25) and of size $4 \times 4$ for Eq. (26). In Eqs. (25) and (26) the letters *B* and *I* denote values of concentration and flux on nodes of the global boundary and the interior interface, while the letter *D* denotes the values of concentration for nodes in the interior of the subregion, since the cell element is discontinuous.

For the solution of the nonlinear problem a linearization procedure is established. In this work, the method of lagging [63] is used. In the initial step of the nonlinear iterative procedure the solution from the previous time step is used. Considering the iterative procedure, the system of equations can be written as,

$$\begin{bmatrix} \mathbf{H}^{(BB)} & \mathbf{H}^{(BI)} \\ \mathbf{H}^{(IB)} & \mathbf{H}^{(II)} \end{bmatrix} \begin{Bmatrix} {}^{k+1}\boldsymbol{\varphi}_{n+1}^{(B)} \\ {}^{k+1}\boldsymbol{\varphi}_{n+1}^{(I)} \end{Bmatrix} = \begin{bmatrix} \mathbf{G}^{(BB)} & \mathbf{G}^{(BI)} \\ \mathbf{G}^{(IB)} & \mathbf{G}^{(II)} \end{bmatrix} \begin{Bmatrix} {}^{k+1}\mathbf{q}_{n+1}^{(B)} \\ {}^{k+1}\mathbf{q}_{n+1}^{(I)} \end{Bmatrix} -$$

$$\frac{1}{\rho}\begin{bmatrix} \mathbf{C}^{(BD)} \\ \mathbf{C}^{(ID)} \end{bmatrix} \left( \frac{1}{\Delta t^a \Gamma(2-a)} \left( \left\{ {}^{k+1}\boldsymbol{\varphi}_{n+1}^{(D)} \right\} - \left\{ \boldsymbol{\varphi}_n^{(D)} \right\} + \left\{ \mathbf{P}_{n,a}^{(D)} \right\} \right) \right. \tag{27}$$

$$\left. - \left\{ \mathbf{LN}\left(\mathbf{m}, {}^{k+1}\boldsymbol{\varphi}_{n+1}^{(D)}\right) \right\} - \left\{ {}^{k+1}\mathbf{f}_{n+1}^{(D)} \right\} \right)$$

$$\left\{ {}^{k+1}\boldsymbol{\varphi}_{n+1}^{(D)} \right\} + \begin{bmatrix} \mathbf{H}^{(DB)} & \mathbf{H}^{(DI)} \end{bmatrix} \begin{Bmatrix} {}^{k+1}\boldsymbol{\varphi}_{n+1}^{(B)} \\ {}^{k+1}\boldsymbol{\varphi}_{n+1}^{(I)} \end{Bmatrix} = \begin{bmatrix} \mathbf{G}^{(DB)} & \mathbf{G}^{(DI)} \end{bmatrix} \begin{Bmatrix} {}^{k+1}\mathbf{q}_{n+1}^{(B)} \\ {}^{k+1}\mathbf{q}_{n+1}^{(I)} \end{Bmatrix} -$$

$$\frac{1}{\rho}\begin{bmatrix} \mathbf{C}^{(DD)} \end{bmatrix} \left( \frac{1}{\Delta t^a \Gamma(2-a)} \left( \left\{ {}^{k+1}\boldsymbol{\varphi}_{n+1}^{(D)} \right\} - \left\{ \boldsymbol{\varphi}_n^{(D)} \right\} + \left\{ \mathbf{P}_{n,a}^{(D)} \right\} \right) \right. \tag{28}$$

$$\left. - \left\{ \mathbf{LN}\left(\mathbf{m}, {}^{k+1}\boldsymbol{\varphi}_{n+1}^{(D)}\right) \right\} - \left\{ {}^{k+1}\mathbf{f}_{n+1}^{(D)} \right\} \right)$$

where $LN\left(\mathbf{m}, {}^{k+1}\varphi_{n+1}\right)$ is the linearized form of Eq. (2) and is given as,

$$LN\left(\mathbf{m}, {}^{k+1}\varphi(\mathbf{r},t)\right) = {}^{k+1}\varphi(\mathbf{r},t)\left(c - d\left({}^{k}\varphi^b(\mathbf{r},t)\right)\right) \tag{29}$$

Consequently, the vector $\mathbf{LN}\left(\mathbf{m}, {}^{k+1}\boldsymbol{\varphi}_{n+1}^{(D)}\right)$ can be expressed as follows,

$$\mathbf{LN}\left(\mathbf{m}, {}^{k+1}\boldsymbol{\varphi}\right) = \left(c[\mathbf{I}] - d\left[\mathbf{M}_1\left({}^{k}\boldsymbol{\varphi}\right)\right]\right)\left\{{}^{k+1}\boldsymbol{\varphi}\right\} = \left[\mathbf{M}\left(\mathbf{m}, {}^{k}\boldsymbol{\varphi}\right)\right]\left\{{}^{k+1}\boldsymbol{\varphi}\right\} \tag{30}$$

Substituting Eq. (30) in Eq. (28) and after some algebra we can obtain the following equation,

$$\left[\mathbf{I} + \mathbf{S}^{(DD)}\right]\left\{{}^{k+1}\boldsymbol{\varphi}_{n+1}^{(D)}\right\} = \begin{bmatrix} \mathbf{G}^{(DB)} & \mathbf{G}^{(DI)} \end{bmatrix} \begin{Bmatrix} {}^{k+1}\mathbf{q}_{n+1}^{(B)} \\ {}^{k+1}\mathbf{q}_{n+1}^{(I)} \end{Bmatrix} -$$
$$\begin{bmatrix} \mathbf{H}^{(DB)} & \mathbf{H}^{(DI)} \end{bmatrix} \begin{Bmatrix} {}^{k+1}\boldsymbol{\varphi}_{n+1}^{(B)} \\ {}^{k+1}\boldsymbol{\varphi}_{n+1}^{(I)} \end{Bmatrix} + \left\{\mathbf{b}^{(DD)}\right\} \tag{31}$$

Where the matrix $\left[\mathbf{S}^{(DD)}\right]$ and the vector $\left\{\mathbf{b}^{(DD)}\right\}$, which change during the nonlinear iteration procedure are given as,

$$\left[\mathbf{S}^{(DD)}\right] = \frac{1}{\rho}\left[\mathbf{C}^{(DD)}\right]\left(c_a[\mathbf{I}] - \left[\mathbf{M}\left(\mathbf{m}, {}^k\boldsymbol{\varphi}_{n+1}^{(D)}\right)\right]\right)$$

$$\left\{\mathbf{b}^{(DD)}\right\} = \frac{c_a}{\rho}\left[\mathbf{C}^{(DD)}\right]\left\{\boldsymbol{\varphi}_n^{(D)}\right\} - \frac{c_a}{\rho}\left[\mathbf{C}^{(DD)}\right]\left\{\mathbf{P}_{n,a}^{(D)}\right\} + \frac{1}{\rho}\left[\mathbf{C}^{(DD)}\right]\left\{{}^{k+1}\mathbf{f}_{n+1}^{(D)}\right\}$$

(32)

To eliminate the volume uknowns using Eq. (31), we can solve for the vector $\left\{{}^{k+1}\boldsymbol{\varphi}_{n+1}^{(D)}\right\}$ and replace the resulting equation in Eq. (27). After some algebra we obtain the following system of equations,

$$\begin{bmatrix} {}^k\overline{\mathbf{H}}^{(BB)} & {}^k\overline{\mathbf{H}}^{(BI)} \\ {}^k\overline{\mathbf{H}}^{(IB)} & {}^k\overline{\mathbf{H}}^{(II)} \end{bmatrix} \begin{Bmatrix} {}^{k+1}\boldsymbol{\varphi}_{n+1}^{(B)} \\ {}^{k+1}\boldsymbol{\varphi}_{n+1}^{(I)} \end{Bmatrix} = \begin{bmatrix} {}^k\overline{\mathbf{G}}^{(BB)} & {}^k\overline{\mathbf{G}}^{(BI)} \\ {}^k\overline{\mathbf{G}}^{(IB)} & {}^k\overline{\mathbf{G}}^{(II)} \end{bmatrix} \begin{Bmatrix} {}^{k+1}\mathbf{q}_{n+1}^{(B)} \\ {}^{k+1}\mathbf{q}_{n+1}^{(I)} \end{Bmatrix} + \begin{Bmatrix} {}^k\overline{\mathbf{b}}_{n+1}^{(B)} \\ {}^k\overline{\mathbf{b}}_{n+1}^{(I)} \end{Bmatrix}$$

(33)

Where the new matrices and the vector introduced have the following form,

$$\begin{bmatrix} {}^k\overline{\mathbf{H}}^{(BB)} & {}^k\overline{\mathbf{H}}^{(BI)} \\ {}^k\overline{\mathbf{H}}^{(IB)} & {}^k\overline{\mathbf{H}}^{(II)} \end{bmatrix} = \begin{bmatrix} \mathbf{H}^{(BB)} & \mathbf{H}^{(BI)} \\ \mathbf{H}^{(IB)} & \mathbf{H}^{(II)} \end{bmatrix} - \begin{bmatrix} \mathbf{S}^{(BD)} \\ \mathbf{S}^{(ID)} \end{bmatrix}\left[\mathbf{I}+\mathbf{S}^{(DD)}\right]^{-1}\begin{bmatrix} \mathbf{H}^{(DB)} & \mathbf{H}^{(DI)} \end{bmatrix}$$

$$\begin{bmatrix} {}^k\overline{\mathbf{G}}^{(BB)} & {}^k\overline{\mathbf{G}}^{(BI)} \\ {}^k\overline{\mathbf{G}}^{(IB)} & {}^k\overline{\mathbf{G}}^{(II)} \end{bmatrix} = \begin{bmatrix} \mathbf{G}^{(BB)} & \mathbf{G}^{(BI)} \\ \mathbf{G}^{(IB)} & \mathbf{G}^{(II)} \end{bmatrix} - \begin{bmatrix} \mathbf{S}^{(BD)} \\ \mathbf{S}^{(ID)} \end{bmatrix}\left[\mathbf{I}+\mathbf{S}^{(DD)}\right]^{-1}\begin{bmatrix} \mathbf{G}^{(DB)} & \mathbf{G}^{(DI)} \end{bmatrix}$$

(34)

$$\begin{Bmatrix} {}^k\overline{\mathbf{b}}_{n+1}^{(B)} \\ {}^k\overline{\mathbf{b}}_{n+1}^{(I)} \end{Bmatrix} = \begin{Bmatrix} \mathbf{b}^{(BD)} \\ \mathbf{b}^{(ID)} \end{Bmatrix} - \begin{bmatrix} \mathbf{S}^{(BD)} \\ \mathbf{S}^{(ID)} \end{bmatrix}\left[\mathbf{I}+\mathbf{S}^{(DD)}\right]^{-1}\left\{\mathbf{b}^{(DD)}\right\}$$

As it is shown, the computation of the new matrices requires the inversion of the matrix $[\mathbf{X}] = \left[\mathbf{I}+\mathbf{S}^{(DD)}\right]$. The inverse of this matrix can be computed iteratively as follows,

$$[\mathbf{X}]_{r+1}^{-1} = [\mathbf{X}]_r^{-1} - g_k [\mathbf{X}]_r^{-1}\left[\mathbf{S}^{(DD)}\right]_r [\mathbf{X}]_r^{-1}, \quad r = 1,\ldots,N_E^{(D)}$$

$$g_k = \frac{1}{1+trace\left(\left[\mathbf{X}_f\right]_r^{-1}\left[\mathbf{S}^{(DD)}\right]_r\right)}$$

(35)

$$[\mathbf{X}]_1^{-1} = [\mathbf{I}]$$

More details about the derivation and the usage of this iterative procedure can be found in Henderson and Searle [64]. Foreach subdomain Eq. (33) is assembled separately and is used to form the final linear system of equation s. After the solution for the boundary nodes of each subdomain is obtained Eq. (31) is used to compute the unknown concentrations $\left\{ {}^{k+1}\boldsymbol{\varphi}_{n+1}^{(D)} \right\}$.

When assembling the final linear system of equations, the interface conditions of Eq. (10) are considered to eliminate the uknown fluxes on the interfaces [56, 61, 65]. Applying also the boundary conditions of the problem we can obtain the following linear system of equations

$$\left[ {}^{k}\mathbf{A} \right] \left\{ {}^{k+1}\mathbf{x}_{n+1} \right\} = \left\{ {}^{k+1}\mathbf{b}_{n+1} \right\} - \left\{ {}^{k}\mathbf{f} \right\} \tag{36}$$

where the matrix $\left[ {}^{k}\mathbf{A} \right]$ appearing in the above equation is sparse and non-symmetric while its sparsity depends on the number of subdomains used for the solution of the problem.

## 4. Numerical results

The accuracy of the proposed method is examined by solving six different problems. When comparing results with a known solution, the following error estimators are used

$$E_\infty = \max_{i=1}^{N} \left| \varphi_i^{exact} - \varphi_i^{num} \right|, \quad E_2 = \frac{\left\| \boldsymbol{\varphi}^{exact} - \boldsymbol{\varphi}^{num} \right\|_2}{\left\| \boldsymbol{\varphi}^{exact} \right\|_2} \tag{37}$$

where the norm appearing in $E_2$ is defined for a vector $\mathbf{v}$ as:

$$\left\| \mathbf{v} \right\|_2 = \left( \sum_{i=1}^{N} v_i^2 \right)^{\frac{1}{2}} \tag{38}$$

In all the following presented problems, units are not mentioned, and it is assumed that all quantities are specified in consistent units. In addition to the provided analytical solutions, for the sake of comparison, an additional solver for the solution of the Fisher-KPP diffusion-reaction equation using the Fragile Points Method (FPM) has been developed. The FPM is a meshless method based on the discontinuous Galerkin weak form and was presented for the first time recently in Dong et al. [57]. Until now it has been used for the solution of various problems in computational mechanics [66-71]. In this work the FPM primal formulation has been implemented [57].

Problem 1

In this problem we consider a closed rectangular domain $\overline{V} = \{x, y : 0 \leq x \leq 1, \ 0 \leq y \leq 2\}$ and compute the solution of Eq. (1) when the fractional time derivative is defined in the Caputo sense. The constant coefficient vector appearing in the nonlinear term is chosen as $\mathbf{m} = [3,1,1]^T$, the diffusivity coefficient is $\rho = 1$, while the function $f(\mathbf{r},t)$ appearing in Eq. (1) for a point $\mathbf{r} = (x, y)$ has the following form,

$$f(\mathbf{r},t) = e^{2x}(1-x^2)t^a \frac{\Gamma(2a+1)}{\Gamma(a+1)} - 2\rho t^{2a}(1-4x-2x^2)e^{2x}$$
$$- \left[t^{2a}(1-x^2)e^{2x}\right]\left[1-\left(t^{2a}(1-x^2)e^{2x}\right)^3\right] \quad (39)$$

The boundary and initial conditions of the problem are the following:

$$\varphi(\mathbf{r},t)\big|_{x=0} = t^{2a}, \; \varphi(\mathbf{r},t)\big|_{x=1} = 0, \; q(\mathbf{r},t)\big|_{y=0} = q(\mathbf{r},t)\big|_{y=2} = 0 \quad (40)$$

$$\varphi(\mathbf{r},0) = 0, \; \mathbf{r} \in V$$

The analytical solution of the problem is provided in Majeed et al. [39]:

$$\varphi(\mathbf{r},t) = t^{2a}(1-x^2)\exp(2x) \quad (41)$$

The domain has been uniformly discretized with 256 quadrilateral elements-subdomains. The solution of the problem has been obtained in the time range $t \in [0, 0.5]$ using a constant time step $\Delta t = 2.5 \times 10^{-3}$. In Fig. 1, the obtained numerical results for the concentration $\varphi$ are presented on the line $y = 0.25$ at $t = 0.50$ considering a constant and a linear approximation of the uknown fields in each subdomain, for the orders of differentiation $a = 0.50 - 0.90$ with a step $\Delta\alpha = 0.10$. It is shown that accurate results can be obtained in both cases. In Table 1 the errors $E_\infty$ and $E_2$, as defined in Eq. (37) are presented for each of the above-mentioned orders of differentiation for the linear mesh results.

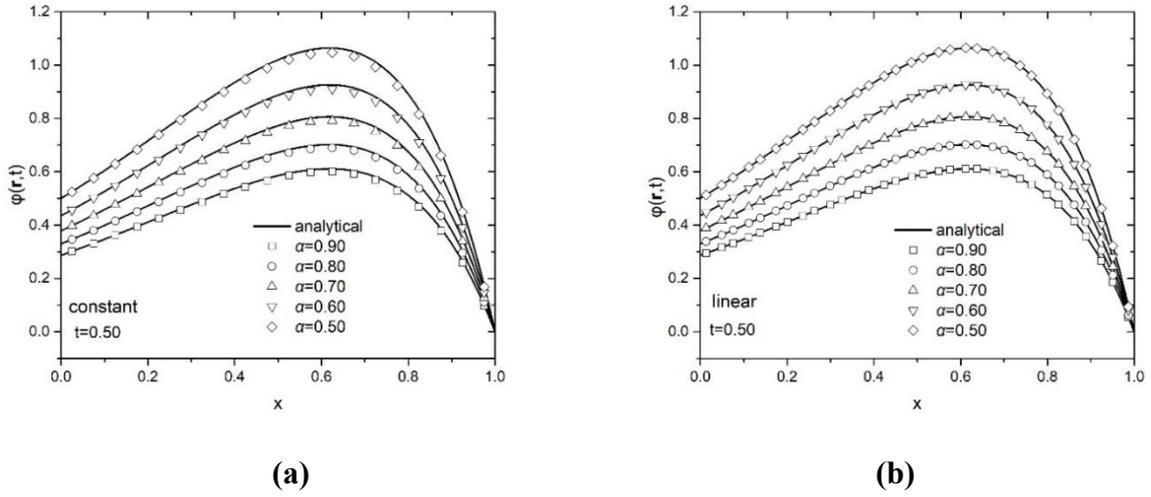

(a)                                    (b)

**Fig. 1** Numerical results for the rectangular domain obtained on the line $y = 0.25$ using, (a) constant elements (b) linear elements.

**Table 1.** Numerical errors on the line $y = 0.25$ for the linear mesh at different times

| $t$ | $E_\infty$ | $E_2$ |
| --- | --- | --- |
| 0.1 | $6 \times 10^{-5}$ | $2.2 \times 10^{-4}$ |
| 0.2 | $1.2 \times 10^{-4}$ | $2.3 \times 10^{-4}$ |
| 0.4 | $2.2 \times 10^{-4}$ | $4.5 \times 10^{-4}$ |
| 0.5 | $2.5 \times 10^{-4}$ | $5.0 \times 10^{-4}$ |

In a next step the convergence of the implemented numerical method is examined by computing the same errors for three different meshes with $N_D = 64,\ 256,\ 1024$ at $t = 0.5$. The results are presented in Table 2, observing that the method is stable and converges, while the accuracy increases with respect to the mesh size.

**Table 2.** Numerical errors along the line $y = 0.25$ at $t = 0.5$ for three different discretizations

| $N_D$ | $E_\infty$ | $E_2$ |
| --- | --- | --- |
| 64 | $3.6 \times 10^{-4}$ | $1.9 \times 10^{-4}$ |
| 256 | $2.5 \times 10^{-4}$ | $1.7 \times 10^{-4}$ |
| 1024 | $2.4 \times 10^{-4}$ | $1.6 \times 10^{-4}$ |

To visualize the different type of solutions that can be obtained, when using a different definition of the fractional derivative, in the following we also examine the solution of Eq. (1) when the time derivative is the fractal fractional derivative in the Riemann-Liouville sense. In Fig. 2 the computed results are presented at $t = 0.5$ for $a = 0.7$ and different values of the $\beta$ parameter. For this case there is not an available analytical solution in the literature, and therefore only the numerical results are presented. From Fig. 2 we note that changes in the $\beta$ parameter leads to nonuniform changes in the concentration function.

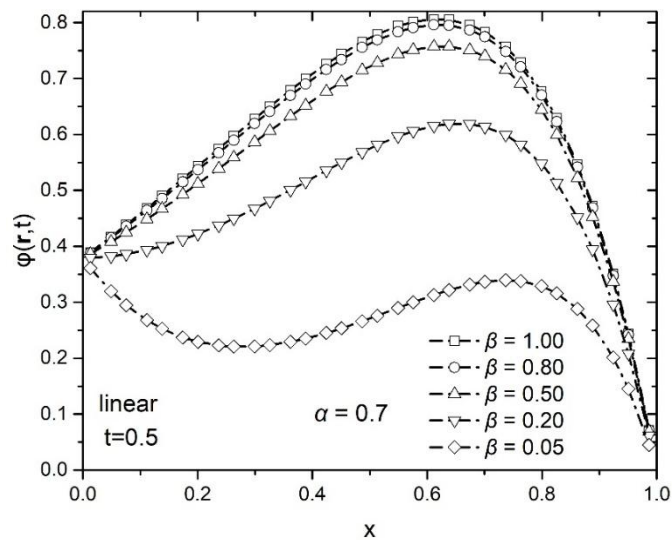

**Fig. 2** Results obtained on the line $y = 0.25$ at $t = 0.5$ using the fractal fractional derivative in the Riemann-Liouville sense.

Problem 2

In this problem we consider initially the case of linear time fractional ($\mathbf{m} = \mathbf{0}$) radial diffusion in a cylinder of radius $R = 2$. For the time derivative, the Caputo fractional derivative is initially considered. The diffusivity value is chosen as $\rho = 1$. The boundary and initial conditions of the problem are the following:

$$\varphi(\mathbf{r},t) = c_0, \quad \mathbf{r} \in \partial V$$

$$\varphi(\mathbf{r},0) = 0, \quad \mathbf{r} \in V$$

(42)

where $c_0$ denotes a constant prescribed value. The analytical solution of the problem is the following [72],

$$\varphi(\mathbf{r},t) = c_0 \left( 1 - 2\sum_{i=1}^{\infty} E_a\left(-\kappa^2 \eta_i^2\right) \frac{J_0(d\eta_i)}{\eta_i J_1(d\eta_i)} \right) \qquad (43)$$

where $d = \dfrac{|\mathbf{r}|}{R}$, $E_a(z) = \sum_{n=0}^{\infty} \dfrac{z^n}{\Gamma(an+1)}$, $a > 0$, $z \in \mathbb{C}$, is the Mittag-Leffler function, $\kappa = \dfrac{\sqrt{\rho}\, t^{a/2}}{R}$ and $n_i = R\xi_i$. The functions $J_0$ and $J_1$ are the Bessel's functions of order zero and one, respectively, while $\xi_i$ are the positive roots of the transcendental equation $J_0(R\xi) = 0$.

The problem has been solved in the time range $t \in [0,1]$ using both the LD-BEM and the FPM, using a time step $\Delta t \approx 3.92 \times 10^{-3}$. Two different discretizations have been prepared, one consisting of 428 quadrilateral elements (Fig. 3), and a second consisting of 900 polygonal elements. The polygonal discretization has been constructed with the Polymesher software [73].

Regarding the implementation of the FPM a primal formulation has been utilized with a weak imposition of the essential boundary conditions [57]. The value of the penalty parameter that has been used corresponds to $n = 5$, and has been chosen in an ad-hoc manner. The integration of the various stiffness terms has been achieved using a one-point integration rule. For the interpolation of the uknown fields the Generalized Finite Difference (GFD) method has been implemented.

In Fig. 4 the normalized concentration is presented as a function of the normalized distance for different orders of the fractional derivative at $t = 1$. These results have been obtained using the mesh of Fig. 3(a). It is shown that the results provided by the two methods are in very good agreement with the analytical solution. In Fig. 5(a) and 5(b) two contour plots are provided corresponding to the results obtained for $a = 0.3$ at $t = 1$. Fig. 5(a) corresponds to the results obtained with the LD-BEM using the quadrilateral mesh, while Fig. 5(b) corresponds to the results obtained with the FPM using polygonal discretization. In Fig. 5(a) the initial quadrilateral mesh has been partitioned in triangles to plot the contour of the result.

In Table 3 we list the numerical errors $E_\infty$ and $E_2$ for both methods for different orders of differentiation at $t = 1$. As it is shown, the two methods are stable and provide accurate results for the entire range of the differentiation order.

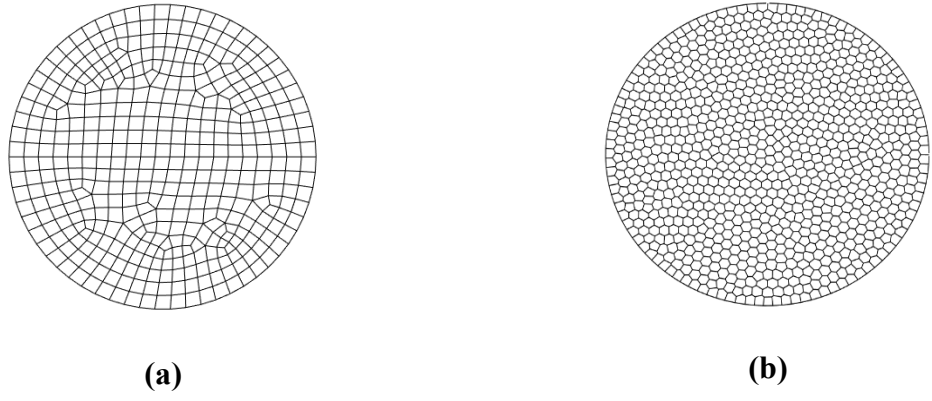

**(a)** **(b)**

**Fig. 3** Discretization of the cylinder with (a) quadrilateral elements and (b) polygonal elements.

**Table 3.** Numerical errors

| | $E_\infty$ | | $E_2$ | |
| --- | --- | --- | --- | --- |
| $a$ | LD-BEM | FPM | LD-BEM | FPM |
| 0.3 | $1.21\times10^{-3}$ | $7.23\times10^{-3}$ | $1.22\times10^{-3}$ | $3.17\times10^{-3}$ |
| 0.5 | $1.47\times10^{-3}$ | $5.36\times10^{-3}$ | $1.37\times10^{-3}$ | $2.27\times10^{-3}$ |
| 0.7 | $1.92\times10^{-3}$ | $3.20\times10^{-3}$ | $1.72\times10^{-3}$ | $1.30\times10^{-3}$ |
| 0.9 | $2.83\times10^{-3}$ | $8.80\times10^{-4}$ | $2.40\times10^{-3}$ | $3.42\times10^{-4}$ |

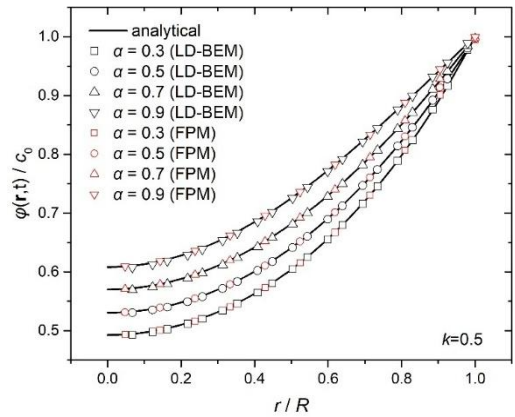

**Fig. 4** Numerical results for both the LD-BEM and the FPM at $t=1$ and various orders of the Caputo fractional derivative.

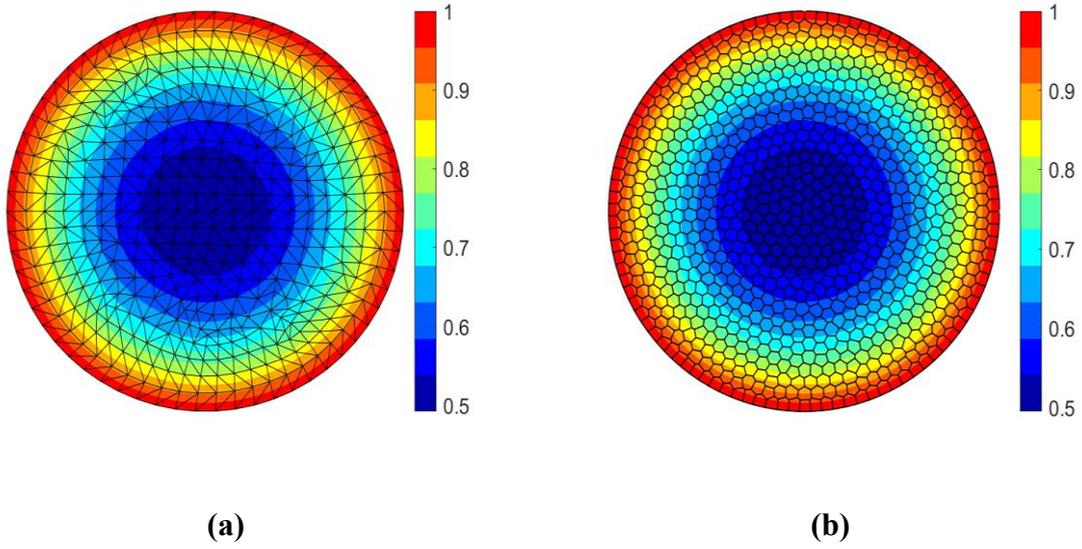

**(a)**  **(b)**

**Fig. 5** Numerical result for the concentration function $\varphi$ at $t=1$ with (a) the LD-BEM with quadrilateral elements and (b) the FPM with polygonal elements.

In Fig. 6 the numerical results are presented for the nonlinear case with $\mathbf{m}=[b,1,1]^T$ and different orders of the polynomial nonlinearity $b$. Here again, the Caputo fractional derivative has been used.

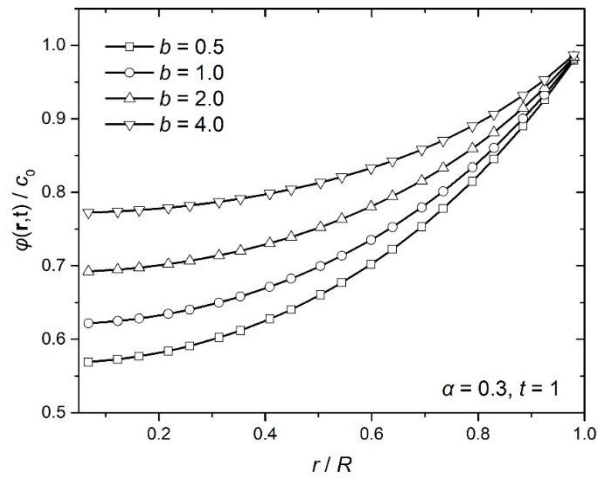

**Fig. 6** Numerical results considering the Caputo fractional derivative with $a=0.3$ and different orders of the polynomial nonlinearity.

Problem 3

In this problem we consider initially the linear ($\mathbf{m}=0$) time-fractional radial diffusion in a hollow cylinder centred at $\mathbf{r}_c=(0,0)$ with internal and external radii $R_{in}=1$ and $R_{out}=2$, respectively. The time fractional derivative is considered in the Caputo sense, while $f(\mathbf{r},t)=0$. The diffusivity is chosen to be $\rho=1$. The boundary and initial conditions of the problem are the following,

$$q(\mathbf{r},t)=0, \ \mathbf{r}\in\partial V_{in}$$

$$\varphi(\mathbf{r},t)=\varphi_0, \ \mathbf{r}\in\partial V_{out} \qquad (44)$$

$$\varphi(\mathbf{r},0)=0, \ \mathbf{r}\in V$$

The analytical solution of the problem in dimensionless form is given as [74]:

$$\varphi_N(\mathbf{r},\tau)=1-\pi\sum_{i=1}^{\infty}\frac{J_1^2(k_i)U_0(dk_i)}{J_1^2(k_i)-J_0^2(\lambda k_i)}E_a\left(-k_i^2\tau^\alpha\right) \qquad (45)$$

where $d=\dfrac{|\mathbf{r}|}{R_{in}}$, $\tau=\dfrac{\rho^{1/\alpha}t}{R_{in}^{2/a}}$, $\lambda=\dfrac{R_{out}}{R_{in}}$, $\varphi_N(\mathbf{r},t)=\dfrac{\varphi(\mathbf{r},t)-\varphi_0}{\varphi_1-\varphi_0}$ are nondimensional quantities and $U_0(dk_i)=J_0(dk_i)Y_0(\lambda k_i)-J_0(\lambda k_i)Y_0(dk_i)$. The functions $J_0$, $Y_0$ and $J_1$, $Y_1$ are Bessel functions of the first and second kind or order zero and one, respectively, while $k_i$ are the positive roots of the transcendental equation $J_1(k_i)Y_0(\lambda k_i)-J_0(\lambda k_i)Y_1(k_i)=0$.

The problem has been solved for $t\in[0,1]$ using both the LD-BEM and the FPM, utilizing a time step $\Delta t=3.92\times 10^{-3}$. Two different spatial discretizations have been prepared, a structured one with 588 quadrilateral elements, and one with 900 polygonal elements. As previously, the polygonal discretization has been constructed with the Polymesher software. The two different discretizations are shown in Fig. 7. For the implementation of the FPM a penalty parameter $n=5$ has been utilized, while for the integration of the various stiffness terms a one-point integration rule has been used. In Fig. 8 the normalized concentration is presented as a function of the normalized distance for different orders of the fractional derivative at $t=1$. In Fig. 9 the results obtained when considering only the linear behavior, are compared with the nonlinear case when using $\mathbf{m}=[1,1,1]^T$ for different orders of differentiation. It is shown that the nonlinear behavior for the chosen parameters leads to similar solution patterns with the linear case. In Fig. 11 the result contours are presented for the nonlinear case with $a=0.3$. The results on the polygonal mesh are obtained with FPM while the results on the quadrilateral mesh are obtained with LD-BEM. The quadrilateral mesh has been partitioned into triangles to plot the results.

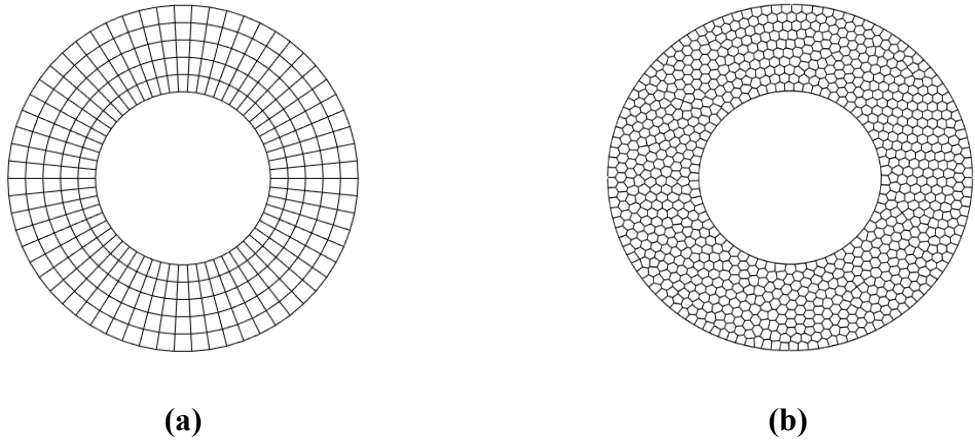

**(a)**                      **(b)**

**Fig. 7** Geometry and discretization with (a) quadrilateral and (b) polygonal elements.

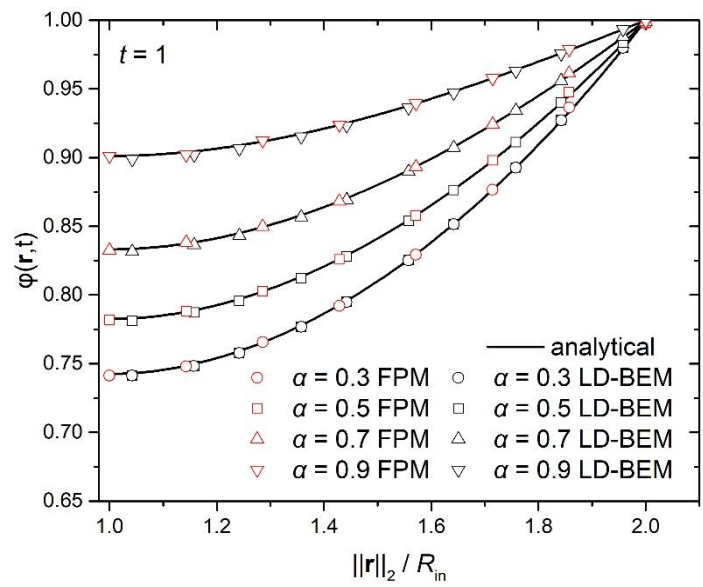

**Fig. 8** Numerical results for the concentration obtained with the LD-BEM and FPM at $t=1$ and different orders of differentiation.

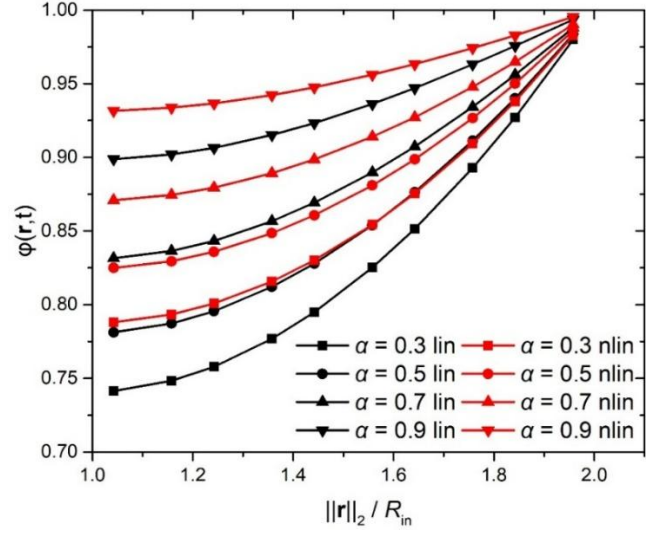

**Fig. 9** Comparison between the linear and the nonlinear results for $\mathbf{m} = [1,1,1]^T$ obtained with LD-BEM considering different orders of differentiation.

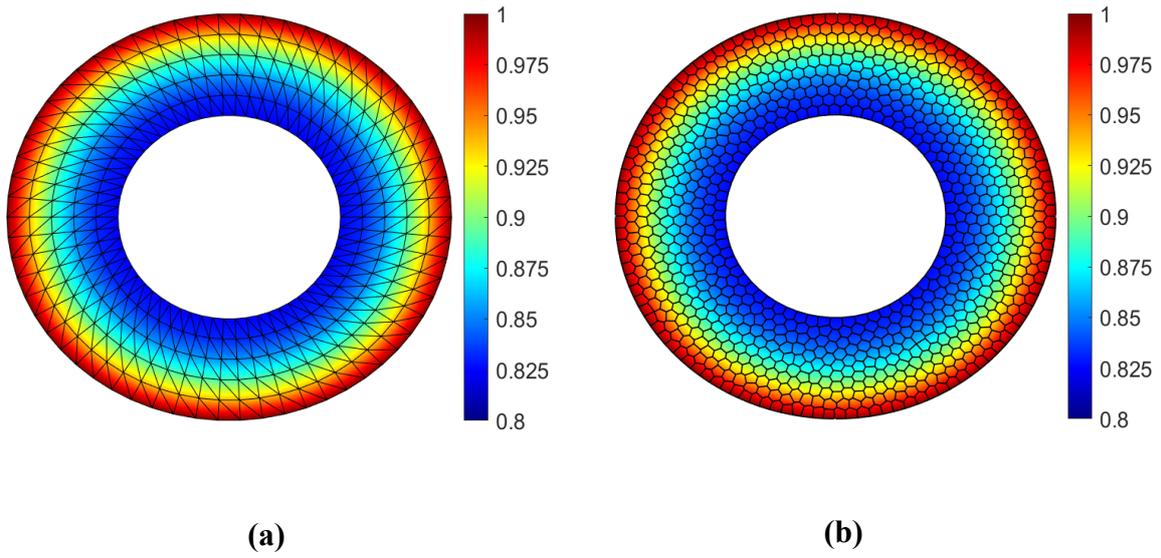

(a)          (b)

**Fig. 10** Numerical results obtained with (a) the LD-BEM and (b) the FPM for $\mathbf{m} = [1,1,1]^T$ and $a = 0.3$.

Problem 4

In this problem we consider the nonlinear diffusion in a rectangular domain $D = \{x, y : -L \leq x \leq L, -B \leq y \leq B\}$ with $L = 10$, $B = 20$. The fractional time derivative appearing in Eq. (1) is defined in the Caputo sense. The constant coefficient vector appearing in the nonlinear partial differential equation is chosen as $\mathbf{m} = [1,1,1]^T$, the diffusivity coefficient is $\rho = 1$, while the function $f(\mathbf{r},t)$ that has been considered has the following form,

$$f(\mathbf{r},t) = \sin x + 2\sin x \frac{t^a}{\Gamma(1+\alpha)} + \sin^2 x \frac{t^{2a}}{(\Gamma(1+a))^2} \tag{46}$$

The boundary and initial condition of the problem are the following,

$$\varphi(\mathbf{r},t)\big|_{x=-L} = 1 + \frac{t^a}{\Gamma(1+\alpha)}\sin(-L), \quad \varphi(\mathbf{r},t)\big|_{x=L} = 1 + \frac{t^a}{\Gamma(1+\alpha)}\sin L$$

$$q(\mathbf{r},t)\big|_{y=-B} = q(\mathbf{r},t)\big|_{y=B} = 0 \tag{47}$$

$$\varphi(\mathbf{r},0) = 1, \quad \mathbf{r} \in V$$

The analytical solution of the problem is the following [34],

$$\varphi(\mathbf{r},t) = 1 + \frac{t^a}{\Gamma(1+\alpha)}\sin x \tag{48}$$

For the discretization of the domain 256 quadrilateral elements and the problem has been solved with LD-BEM in the time range $t \in [0, 0.1]$, considering a time step $\Delta t = 1.01 \times 10^{-3}$. In Fig. 11 the numerical results for the concentration are presented along the line $y = 10$ at $t = 0.1$ for different orders of differentiation. It is shown the results are in good agreement with the analytical solution.

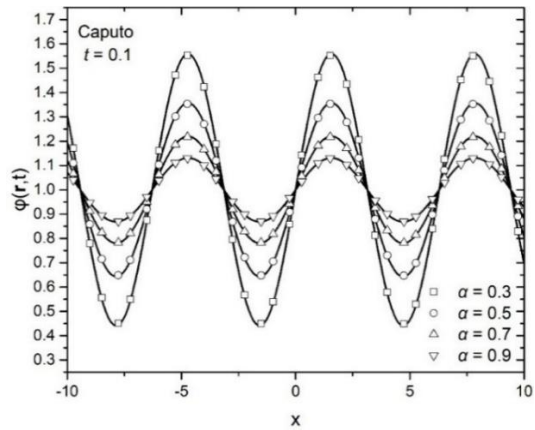

**Fig. 11** Numerical results along the line $y = 10$ at $t = 1$ and various orders of differentiation.

Problem 5

In this problem we consider the nonlinear time fractional diffusion in the domain shown in Fig. 12. This domain is constructed as the intersection of two circles, where the first is centred in $x_{c_1}=(0,0)$ with a radius $R_1=1$ while the second circle is centred in $x_{c_2}=(-0.4,0)$ with a radius $R_2=0.55$. The fractional time derivative appearing in Eq. (1) is defined in the Caputo sense. The constant coefficient vector appearing in the nonlinear partial differential equation is chosen as $\mathbf{m}=[2,1,1]^T$, while the diffusivity coefficient is set to $\rho=1$. The boundary and the initial conditions of the problem are the following

$$\varphi(\mathbf{r},t)\big|_{S_{R_1}}=1,\ q(\mathbf{r},t)\big|_{S-S_{R_1}}=0$$

$$\varphi(\mathbf{r},0)=0,\ \mathbf{r}\in V$$

(49)

where $S_{R_1}$ is shown in Fig. 12. Two different domain discretizations have been constructed. The first consists of 880 quadrilateral elements and the second consists of 900 polygonal elements. The problem has been solved in the time range $t\in[0,0.2]$ while a time step of $\Delta t=7.84\times10^{-4}$ has been considered. In Fig. 13 the concentration results are presented along the line $\{y=0,x\geq0\}$ at $t=0.2$ for three different orders of differentiation. It is shown that the agreement of the results obtained with the two methods is excellent. In Fig. 14 two different contour plots are provided with the results obtained with the two methods for $a=0.3$. In Fig. 15 the numerical results are presented for different orders of the polynomial nonlinearity considering a fixed order of differentiation $a=0.5$.

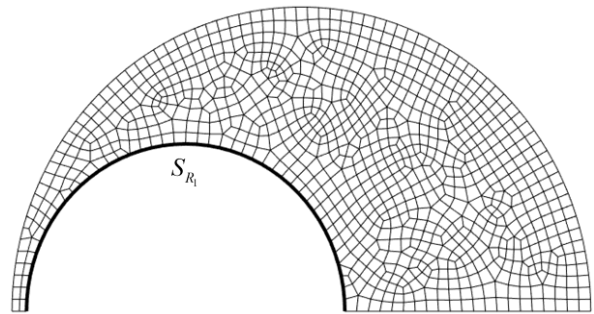

**Fig. 12** Geometry of the problem and discretization with quadrilateral elements.

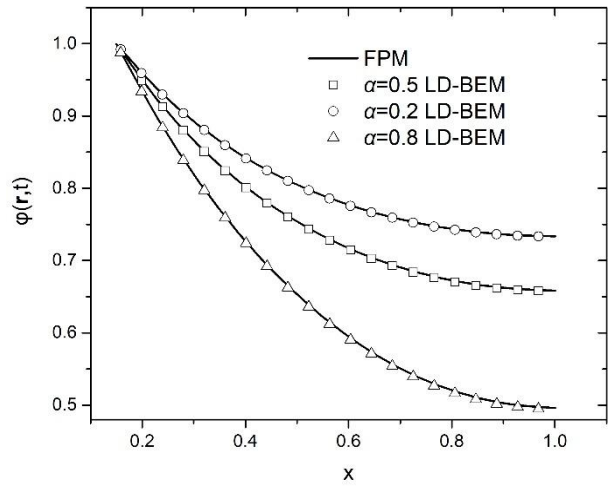

**Fig. 13** Numerical solution for the set of nodes $\{y = 0, x > 0\}$ for different values of the differentiation order.

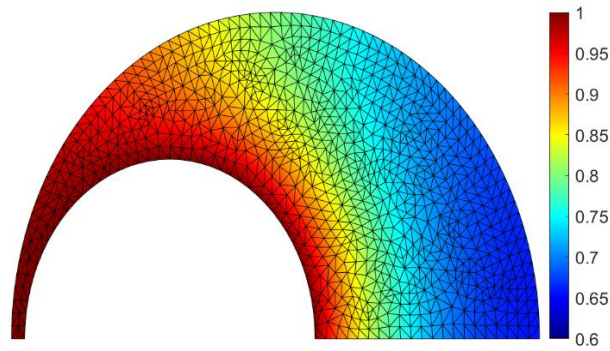

(**a**)

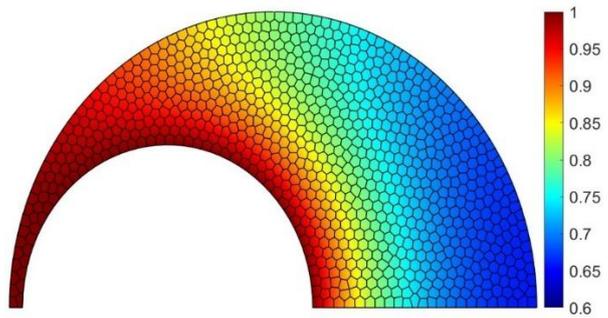

(**b**)

**Fig. 14** Contour plots for the solution obtained using (a) the LD-BEM and (b) the FPM for $a = 0.3$.

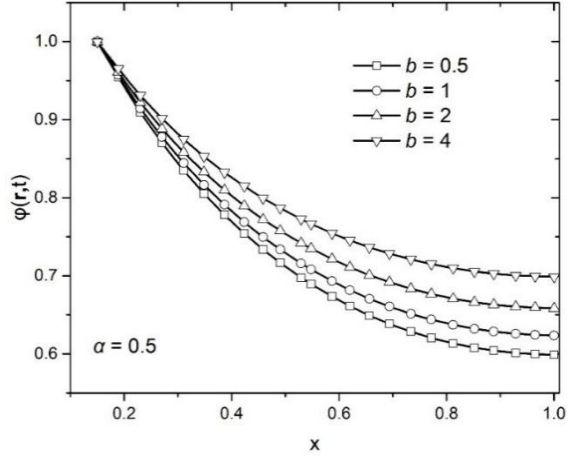

**Fig. 15** Numerical solution considering different orders of the polynomial nonlinearity.

Problem 6

In this problem we consider the time fractional linear and nonlinear diffusion on the domain shown in Fig. 16. The radius for the points $\left[r(\theta),\theta\right]$ on the boundary of this geometry is given as:

$$r(\theta) = R + 0.2R\cos(3\theta) + 0.02R\cos(5\theta) + \\ 0.4R\sin(8\theta) + 0.4R\sin(4\theta) + 0.1R\sin(15\theta) \quad (50)$$

with $R = 0.4$. The geometry of the domain has been generated using 200 points with $\theta \in [0, 2\pi)$ connected with straight lines and has been discretized with a mesh consisting of 2181 quadrilateral elements. The boundary and the initial conditions of the problem are given as follows

$$\varphi(\mathbf{r},t)\big|_{S_{R_1}} = 10, \quad q(\mathbf{r},t)\big|_{S-S_{R_1}} = 0 \quad (51)$$

$$\varphi(\mathbf{r},0) = 0, \quad \mathbf{r} \in V$$

The problem has been solved with $t \in [0,10]$ and a time step $\Delta t = 3.922 \times 10^{-2}$. The diffusivity is chosen to be $\rho = 1$. As in the previous examples, in the FPM setup a penalty parameter $n = 5$ has been used. In Fig 17 the numerical results obtained with the two methods, considering the case of linear time fractional diffusion ($\mathbf{m} = \mathbf{0}$), with different orders of differentiation, are plotted along the line $L$ shown in Fig. 16 at $t = 10$. In Fig. 18 the results for the concentration are presented for the nonlinear time fractional diffusion with $\mathbf{m} = [b,1,1]^T$ considering three different orders of the polynomial nonlinearity, while In Fig. 19 and in Fig 20 the contour plots are provided at $t=10$ for the case of linear time fractional diffusion with $a=0.3$ and for the case of non-linear time fractional diffusion with the same differentiation order and for the constant vector $\mathbf{m}=\begin{bmatrix}1,1,1\end{bmatrix}^T$.

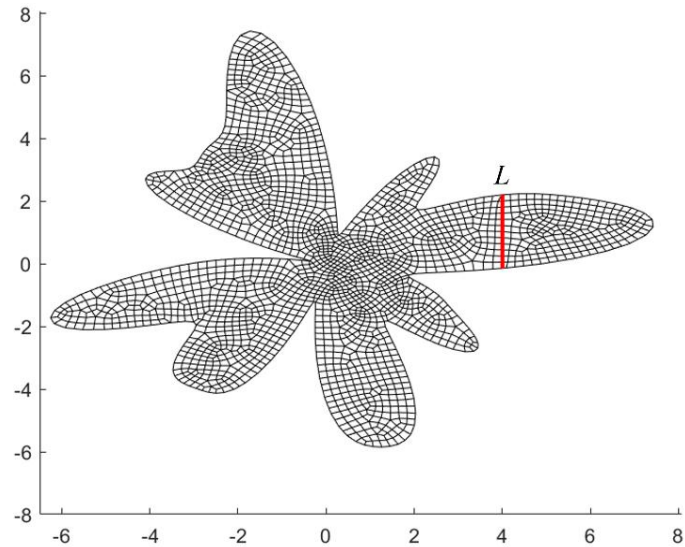

**Fig. 16** Geometry and discretization of the problem

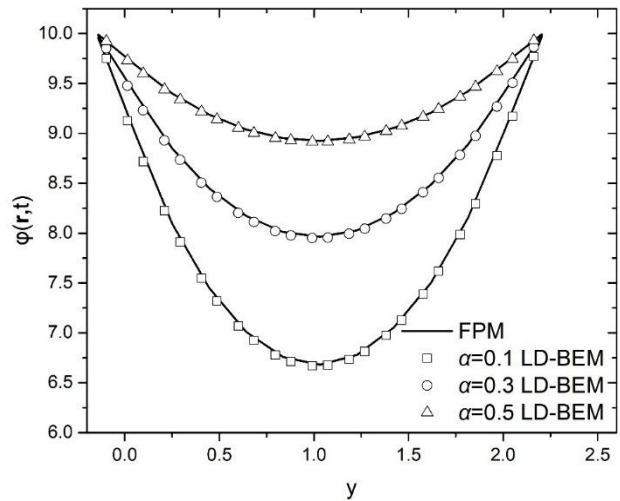

**Fig. 17** Numerical solution obtained with LD-BEM and FPM for $\mathbf{m}=\mathbf{0}$ at $t=10$ on the line $L$ for different orders of differentiation

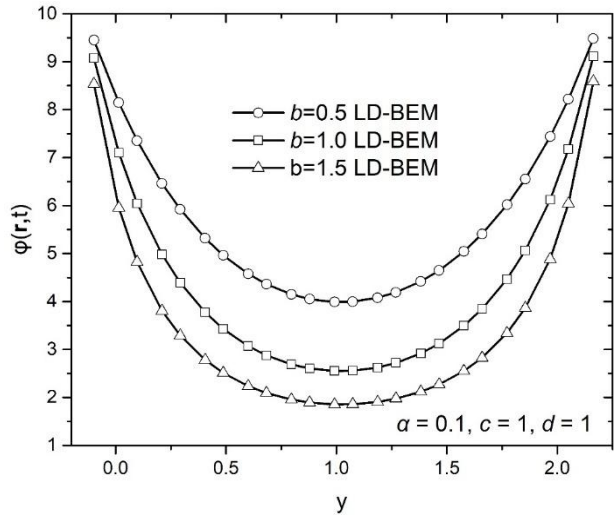

**Fig. 18** Numerical solution obtained with $\mathbf{m} = [b,1,1]^T$ for different orders of the polynomial nonlinearity

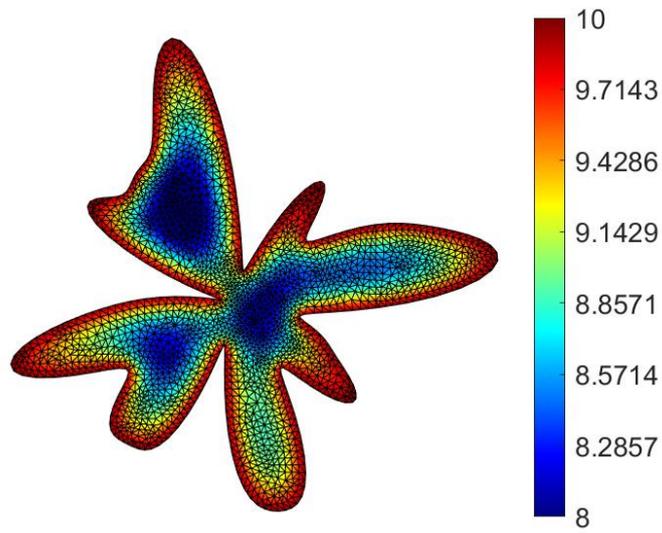

**Fig. 19** Distribution of the concentration $\varphi$ for the case $\mathbf{m} = \mathbf{0}$ at $t = 10$.

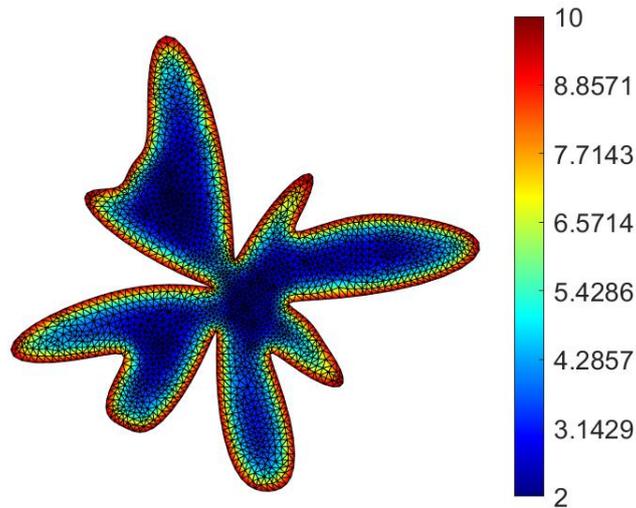

**Fig. 20** Distribution of the concentration $\varphi$ for the case $\mathbf{m} = [1,1,1]^T$ at $t = 10$.

## 5. Conclusions

In this work, an extension of the local domain boundary element method for solving the nonlinear time fractional Fisher-KPP equation has been presented. Solutions for a variety of two-dimensional problems with domains of varying complexity and different orders of the polynomial nonlinearity have been presented using the Caputo fractional derivative and a fractal fractional derivative in the Riemann-Liouville sense. The accuracy of the method has been verified through analytical solutions while comparisons have also been performed by implementing a version of the recently proposed FPM method. It has also been shown that terms such as volume sources can be easily considered in the formulation of the LD-BEM. Extensions to other type of fractional derivatives can be implemented in a straightforward manner.

**Acknowledgement**

The research work has been supported by the Hellenic Foundation for Research and Innovation (HFRI) under the "First Call for HFRI. Research Projects to support Faculty members and Researchers and the procurement of high-cost research equipment grant" (Project Number: 2060).